\documentclass{article}

\usepackage{amsfonts}
\usepackage{graphics}
\usepackage{amssymb}
\usepackage{amsfonts,amsmath}
\usepackage[all]{xy}

\newtheorem{theorem}{Theorem}[section]
\newtheorem{lemma}{Lemma}[section]
\newtheorem{definition}{Definition}[section]
\newtheorem{proposition}[theorem]{Proposition}
\newtheorem{corollary}[theorem]{Corollary}

\newcommand{\eop }{ \hfill $\Box$ }

\begin{document}
\begin{center}
{\Large An equivalence between harmonic sections and sections that are harmonic maps\\}

\end{center}

\vspace{0.3cm}

\begin{center}
{\large  Sim\~ao Stelmastchuk}  \\

\textit{Departamento de  Matem\'{a}tica, Universidade Estadual de Campinas,\\ 13.081-970 -  Campinas - SP, Brazil. e-mail: simnaos@gmail.com}
\end{center}

\vspace{0.3cm}

\begin{abstract}
Let $\pi:(E,\nabla^{E}) \rightarrow (M,g)$ be an affine submersion with horizontal distribution, where $\nabla^{E}$ is a symmetric connection and $M$ is a Riemannian manifold. Let $\sigma$ be a section of $\pi$, namely, $\pi \circ \sigma = Id_{M}$. It is possible to study the harmonic property of section $\sigma$ in two ways. First, we see $\sigma$ as a harmonic map. Second, we see $\sigma$ as harmonic section. In the Riemannian context, it means that $\sigma$ is a critical point of the vertical functional energy. Our main goal is to find conditions to the assertion: $\sigma$ is a harmonic map if and only if $\sigma$ is a harmonic section.
\end{abstract}


\vspace{0.3cm} \noindent {\bf MSC2010 subject classification:}
 53C43, 58E20, 58J65, 60H30.

\section{Introduction}
Let $M$ be a Riemannian manifold and $E$ a differential manifold with a symmetric connection $\nabla^{E}$. Let us denote the Levi-Civita connection on $M$ by $\nabla^{M}$. Let $\pi:E \rightarrow M$ be a submersion. Then we can define the vertical bundle in $TE$ as $VE=\rm{Ker}\,(\pi_{*})$. Let $HE$ be a smooth distribution in $TE$ such that $TE= VE \oplus HE$. We call $HE$ of horizontal distribution. Let us denote by $\mathbf{v}:TE \rightarrow VE$ and $\mathbf{h}:TE \rightarrow HE$ the  vertical and horizontal projectors, respectively. We denote by $H_{x}=\pi_{p}|_{H_{p}}^{-1}:T_{x}M \rightarrow H_{p}E$, where $\pi(p)=x$, the family of isomorphism which we call the horizontal lift. The submersion $\pi$ is called affine submersion with horizontal distribution if $\mathbf{h}\nabla^{E}_{H(X)}HY = H(\nabla^{M}_{X}Y)$ for the $X,Y$ vector fields on $M$ (see \cite{abe} for more details).
 
Let $\sigma$ be a section of an affine submersion with horizontal distribution $\pi$, that is, $\pi \circ \sigma = Id_{M}$. We can study the harmonicity of $\sigma$ in two ways. First way, $\sigma$ is a harmonic map. Let $\{e_{1},\ldots,e_{n}\}$ be an arbitrary local orthonormal frame fields at $x \in M$. The tension field of $\sigma$ is given by
\[
 \tau_{\sigma}(x) = \sum_{i=1}^{n} (\nabla^{E}_{\sigma_{*}e_{i}}\sigma_{*}e_{i} - \sigma_{*}\nabla^{M}_{e_{i}}e_{i})(x).
\]
We say that $\sigma$ is a harmonic map if $\tau_{\sigma} \equiv 0$. We observe that $\nabla^{E}$ is no necessarily the Levi-Civita connection. 

The second way is to study the harmonic property only through the vertical component of $\sigma$. Let $\nabla^{v}$ be a vertical connection on $E$, namely, $\nabla_{X}V= \mathbf{v}\nabla^{E}_{X}V$, where $X$ is a vector field and $V$ is a vertical vector field on $E$.  Let $\mathbf{v}\sigma_{*}$ be the vertical component of $\sigma_{*}$. The vertical tension field is given by
\[
 \tau_{\sigma}^{v}(x) = \sum_{i=1}^{n} (\nabla^{v}_{\mathbf{v}\sigma_{*}e_{i}}\mathbf{v}\sigma_{*}e_{i} - \mathbf{v}\sigma_{*}\nabla^{M}_{e_{i}}e_{i})(x).
\]
We say that $\sigma$ is a harmonic section if $\tau_{\sigma}^{v} \equiv 0$. This definition is an extension of one gived by C. Woods \cite{wood1} when $\pi$ is a Riemannian submersion.

Our main objective is to find conditions for the assertion:
\begin{center}
$\sigma$ is a harmonic map if and only if $\sigma$ is a harmonic section.
\end{center}

At first, seen by a geometric point of view it is not always true (see proposition \ref{hsprop2}), unless we stablish some conditions about $\sigma$. Our wish is to find geometric conditions that satify the assertion. With this aim, we use the Fundamental tensors $A$ and $T$ given by B. O'Neill \cite{oneill}(see condition (\ref{skewsymmetry})). Under some conditions in $A$ and $T$ it is possible to prove with stochastic calculus  the result below (see Theorem \ref{hsteo1}).\\

{\bf Theorem:}
 $\sigma$ is a harmonic section if and only if $\mathbf{v}\tau_{\sigma} \equiv 0$.\\

This result is important because it gives a good condition to the assertion. From this and Proposition \ref{hsprop2} follows our principal result:\\

{\bf Theorem :}
{\it
Let $\sigma$ be a section of $\pi$ and $C_{\sigma}$ and $D_{\sigma}$ be tensors defined by (\ref{hsprop2eq2})}. 
\begin{enumerate}
\item If $\rm{tr}\, D_{\sigma}= 0$, then $\sigma$ is a harmonic map if and only if $\sigma$ is a harmonic section;
\item If $\sigma$ is a harmonic map, then $\rm{tr}\, C_{\sigma} = 0 $ and $\rm{tr}\, D_{\sigma}=0$;
\item If $\rm{tr}\, C_{\sigma} \neq 0$ or $\rm{tr}\, D_{\sigma} \neq 0$,  then $\sigma$ is not a harmonic map;
\item If $\sigma$ is a harmonic section, then $\rm{tr}\, C_{\sigma} = 0$;
\item If $\rm{tr}\, C_{\sigma} \neq 0$, then $\sigma$ is not a harmonic section.\\
\end{enumerate}

Item (1) is the answer for our assertion. The other items are good conditions to find candidates that are harmonic sections or sections which are candidate to be harmonic maps. They can be thought as critical points in the sections of $\pi$ set.

We apply this result in many contexts: product manifolds, tangent bundle with horizontal lift, submersion in the sense of Blumenthal and Riemannian submersion. Latter, we show that the geometric condition to the assertion is that the horizontal distribution given by $HE$ is integrable (see Proposition \ref{exprop4}).

%
%
\section{Preliminaries}

We begin by recalling some fundamental facts on Schwartz geometry and stochastic calculus on manifolds. We shall use freely concepts and notations from S. Kobayashi and N. Nomizu \cite{kobay}, L. Schwartz \cite{schwartz}, P.A. Meyer \cite{meyer} and M. Emery \cite{emery1}. A quick survey in these concepts is described by P. Catuogno in \cite{catuogno2}.

Let $M$ be a smooth manifold and $x \in M$. The second order tangent space to $M$ at $x$, which is denoted by $\tau_{x}M$, is the vector space of all differential operators on M at x of order at most two without a constant term. Let $(x_{1}, \ldots, x_{n})$ be a local system of coordinates. Every $L \in \tau_{x}M$ can be written in a unique way as 
\[
 L = a_{ij}D_{ij} + a_{i}D_{i},
\]
where $a_{ij}=a_{ji}$, $D_{i}= \frac{\partial}{\partial x^{i}}$ and $D_{ij}=\frac{\partial^{2}}{\partial x^{i}\partial x^{j}}$ are differential operators at x (we shall use the convention of summing over repeated indices). The elements of $\tau_{x}M$ are called second order tangent vectors at $x$, the elements of the dual vector space $\tau_{x}^{*}M$ are called second order forms at $x$.

The disjoint union $\tau M = \bigcup_{x \in M} \tau_{x}M$ (respectively, $\tau^{*}M  = \bigcup_{x \in M} \tau_{x}^{*}M$) is canonically endowed with a vector bundle structure over $M$, which is called the second order tangent fiber bundle (respectively, second order cotangent fiber bundle) of $M$.

Let $\pi:E \rightarrow M$ be a submersion. Suppose that $TE = VE \oplus HE$, where $VE= \rm{ker}\,(\pi_{*})$ is the vertical distribution and $HE$ is a horizontal distribution of $TM$, respectively. A 1-form $\theta$ in $T^{*}E$ is called vertical form if $\theta(X) = 0$ when $X \in HE$. Let $\tau E = V\tau E \oplus W$, where $V\tau E = \rm{Ker}\,(\pi_{*})$ is the vertical distribution and $W$ is the complementar bundle of $V\tau E$ in $\tau E$. A second order form $\Theta$ on $M$ is also called vertical if the restriction of $\Theta$ to $W$ is null.

Let $M, N$ be smooth manifolds, $F :M\rightarrow N$ a smooth map and $L\in \tau _{x}M$. The differential of $F$ , $F_{*}(x):\tau_{x}M \rightarrow \tau_{F(x)}N$, is given by 
\[
F _{*}(x)L(f)=L_{x}(f\circ F ),
\]
where $f \in \mathcal{C}^{\infty}(N)$.


Let $L$ be a second order vector field on $M$. The square operator of $L$, denoted by $QL$, is the symmetric tensor given by 
\[
QL(f,g)=\frac{1}{2}(L(fg)-fL(g)-gL(f)),
\]
where $f,g\in C^{\infty }(M)$. Let $x \in M$. We consider $Q_{x}:\tau_{x}M\rightarrow T_{x}M\odot T_{x}M$ as the linear application defined by
\[
Q_{x}(L = a_{ij}D_{ij} + a_{i}D_{i}) = a_{ij}D_{i} \odot D_{j}. 
\]

Push forward of second order vectors by smooth maps is related to the so called Schwartz morphisms
between second order tangent vector bundles.

\begin{definition}
Let $M$ and $N$ be smooth manifolds, $x\in M$ and $y\in N$. A linear application $F:\tau_{x}M\rightarrow \tau _{y}N$ is called Schwartz morphism if 
\begin{enumerate}
\item  $F(T_{x}M)\subset T_{y}N$;
\item  for all $L\in \tau _{x}M$ we have $Q(FL)=(F\otimes F)(QL)$.
\end{enumerate}
\end{definition}

A linear application $F:\tau_{x}M \rightarrow \tau_{y}N$ is a Schwartz morphism if and only if there exists a smooth map $\phi: M \rightarrow N$ with $\phi(x)=y$ such that $F=\phi_{x*}$ (see for example Proposition 1 in \cite{emery2}).

Let $(\Omega, (\mathcal{F}_t),\mathbb{P})$ be a filtered probability space which satisfies the usual conditions (see for instance \cite{emery1}).
\begin{definition}
Let $M$ be a smooth manifold. Let $X$ be a stochastic process with values in $M$. We call $X$ a semimartingale if, for all $f$ smooth on $M$, $f(X)$ is a real semimartingale.
\end{definition}

L. Schwartz has noticed, in \cite{schwartz}, that, if $X$ is a continuous semimartingale in a smooth manifold $M$, the It\^o's differentials $dX_{i}$ and $d[X_{i},X_{j}]$ (where $(x_{i})$ is a local coordinate system and $X_{i}$ is the ith coordinate of $X$ in this system) behave under a change of coordinates as the coefficients of a second order tangent vector. The (purely formal) stochastic differential
\[
 d^{2}X_{t} = dX^{i}_{t}D_{i} + \frac{1}{2}d[X^{i},X^{j}]_{t}D_{ij},
\]
is a linear differential operator on $M$, at $X_{t}$, of order at most two, with no constant term. Therefore, the tangent object to $X_{t}$ is formally one of second order. This fact is known as Schwartz principle.

From now on we assume that all semimartingales are continuos.

Let $A$ and $B$ be vector fields on $M$. We call $\{A, B\}= \frac{1}{2}(AB + BA)$ the commutator of the vector fields. We observe that $\{A, B\}$ is a second order vector field. The set of commutators in $\tau M$ will be denoted by $\mathrm{Com}(M)$.
\begin{lemma}\label{plem1}
$Q|_{\mathrm{Com}(M)}: \mathrm{Com}(M) \rightarrow TM \odot TM$ is a one-to-one correspondence.
\end{lemma}
\begin{proof}
It is sufficient  to define the map $\bar{Q}: TM \odot TM \rightarrow \tau M$ by
\[
\bar{Q}(A \odot B) = \{A,B\}.
\]
It is immediate that $Q\left|_{\mathrm{Com}(M)}\right.$ and $\bar{Q}$ are inverse maps.
\eop
\end{proof}

\begin{proposition}\label{pprop1}
Let $M, N$ be smooth manifolds and $F:M \rightarrow N$ be a smooth map. Then $F_{*}:\tau M \rightarrow \tau N$ is written as 
\[
F_{*} = F_{*} \oplus F_{*}\odot F_{*}.
\]
\end{proposition}
\begin{proof}
Let $X$ be a semimartingale in $M$ and $(x_{1}, \ldots, x_{n})$ be a local coordinate system. By Schwartz principle, it follows that 
\[
d^{2}X_{t}= dX^{i}D_{i} + d[X_{i},X_{j}]\{D_{i},D_{j}\}.
\]
Applying $F_{*}$ in $d^{2}X_{t}$ it yields that
\begin{eqnarray*}
F_{*}(d^{2}X_{t}) & = & dX^{i}F_{*}(D_{i}) + d[X_{i},X_{j}]F_{*}(\{D_{i},D_{j}\}).
\end{eqnarray*}
Since $F_{*}$ is a Schwartz morphism, applying $Q$ in $F_{*}(d^{2}X_{t})$ we obtain
\[
Q(F_{*} d^{2}X_{t}) = d[X^{i},X^{j}]F_{*}\otimes F_{*}(D_{i}\odot D_{j}).
\]
Therefore
\[
F_{*}(d^{2}X_{t})  =  dX^{i}F_{*}(D_{i}) + \frac{1}{2}d[X^{i},X^{j}]F_{*}\otimes F_{*}(D_{i}\odot D_{j}),
\]
and the Proposition follows.
\eop
\end{proof}


Let $X$ be a semimartingale in $M$. Let $\Theta_{X_t} \in \tau^*_{X_t}M$ be an adapted stochastic second order form along $X_t$. Let $(U,x^i)$ be a local coordinate system in $M$. With respect to this chart the second order form $\Theta$ can be written as $\Theta_x= \theta_i(x)d^2x^i + \theta_{ij}(x)dx^i\cdot dx^j$ where $\theta_i$ and $\theta_{ij}=\theta_{ji}$ are ($\mathcal{C}^{\infty}$ say) functions in $M$. Then the integral of $\Theta$ along $X$ is defined by:
\[
\int_0^t \Theta ~d^2X=\int_0^t \theta_i(X_s)dX^i_s + \int_0^t \theta_{ij}(X_s)d[X^i,X^j]_s
\]

Let $b$ be a section of $T^{2}_{0}(M)$, which is defined along $X$. The quadratic integral of $b$ along $X$ is defined, locally, by
\[
 \int_0^{t}b\;(dX,dX) = \int_{0}^{t} b_{ij}(X_{t}) d[X^{i},X^{j}]_{t}, 
\]
where $b(x) = b_{ij}(x)dx^{i} \otimes dx^{j}$ and $b_{ij}$ are smooth functions.

Let $M$ be a smooth manifolf endowed with symmetric connection $\nabla^{M}$. In \cite{meyer},  P. Meyer showed that for $\nabla^{M}$ there exists a section $\Gamma^{M}$ in $Hom(\tau M,TM)$ such that $\Gamma^{M}|_{TM}=Id_{TM}$ and  $\Gamma^{M}(AB)= \nabla^{M}_{A}B$, where $A,B \in TM$. We also call $\Gamma^{M}$ by connection.

P. Catuogno gived the following definition in \cite{catuogno2}.

\begin{definition}\label{formafundamentaleharmonica}
Let $M$ and $N$ be smooth manifolds endowed with connections $\Gamma^M$ and $\Gamma^N$, respectively, and $F :M\rightarrow N$ be a smooth map. We define the section $\alpha_F$ in  $\tau^*M\otimes F^*TN$ by 
\begin{equation}\label{fundamental}
\alpha_F=\Gamma^N\circ F_*-F_* \circ \Gamma^M.
\end{equation}
The second fundamental form of $F$, denoted by $\beta_F$, is the unique section of $ (TM \odot TM)^*\otimes F^*TN$ such that $\alpha_F = \beta_F\circ Q$. 

In the case that $(M,g)$ is a Riemmanian manifold and $\Gamma^M$ is the Levi-Civita connection, the tension field of $F$, $\tau_F:M \rightarrow TN$, is given by 
\[
\tau_F=tr\beta_F.
\]
The map $F$ is called harmonic with respect to $\Gamma^{N}$ if its tension field is null, i.e., $\tau_{F}=0$.
\end{definition}

Let $M$ be a smooth manifold endowed with a symmetric connection $\Gamma^{M}$. Let $X$ be a semimartingale in $M$ and $\theta$ be a 1-form along $X$. The It\^o integral is defined by $\int_{0}^{t} \theta d^{M}X:=\int_{0}^{t}\Gamma^{M*}\theta d^{2}X_{t}$.

\begin{definition}
Let $M$ be a Riemmanian manifold with metric $g$. A semimartingale $B$ in $M$ is called Brownian motion if $\int_{0}^{t} \theta d^{M}B_t$ is a real local martingale for all $\theta \in T^{*}M$, where $\Gamma^{M}$ is the Levi-Civita connection, and for any section $b$ in $T^{2}_{0}(M)$ we have 
\begin{equation}\label{Brownian}
\int_0^tb(dB,dB)=\int_0^t\rm{tr}\,b_{B_{s}}ds.
\end{equation}
\end{definition}

%

%
%

\section{An equivalence between harmonic sections and sections that are harmonic maps}

Let $E$ be a diffenrential manifold and $M$ be a Riemannian manifold such that there is a smooth submersion $\pi:E \rightarrow M$. Let $\nabla^{E}$ be a symmetric connection and $\nabla^{M}$ be the Levi-Civita connection on $E$. Let $VE=\mathrm{ker}\,(\pi_{*})$ be the vertical distribution and $HE$ a smooth distribution in $TE$ such that $TE= VE \oplus HE$. Let $\mathbf{v}: TE \rightarrow VE$ and $\mathbf{h}:TE \rightarrow HE$ be the vertical and horizontal projectors, respectively. Let $H_{x} = (\pi_{p∗}|_{H_{p}E})^{-1} : T_{x}M \rightarrow H_{p}E$ be the horizontal lift, where $\pi(p)=x$. We observe that $H_{x}$ is an isomorphism for each $x \in M$. The submersion $\pi:E \rightarrow M$ is called affine submersion with horizontal distribution if  $\mathbf{h}\nabla^{E}_{H(X)}H(Y) = H(\nabla^{M}_{X}Y)$ (see \cite{abe} for more details). A Riemmanian submersion is a classical example of affine submersion with horizontal distribution. Unlles otherwise stated, we assume that $\pi:E \rightarrow M$ is an affine submersion with horizontal distribution. 

We follow the definition of B. O'Neill \cite{oneill} for the Fundamental tensors $A$ and $T$. They are defined by
\[
 T_{U}V = \mathbf{h} \nabla^{E}_{\mathbf{v}{U}}\mathbf{v}{V} + \mathbf{v}\nabla^{E}_{\mathbf{v}U} \mathbf{h}V
\]
and
\[
 A_{U}V = \mathbf{v} \nabla^{E}_{\mathbf{h}{U}}\mathbf{h}{V} + \mathbf{h}\nabla^{E}_{\mathbf{h}U} \mathbf{v}V,
\]
where $U,V$ are vector fields on $E$.

The following assumption will be needed throughout the paper. Let $X,Y$ be vector fields on $M$. We ask connections $\nabla^{E}$ on $E$ such that 
\begin{equation}\label{skewsymmetry}
A_{H(X)}H(Y) = - A_{H(Y)}H(X). 
\end{equation}
We observe that Riemannian submersions satisfy this condition.

The next Lemma is important in the sequence.
\begin{lemma}
Let $\pi:E \rightarrow M$ be a submmersion. Let $(x_{1}, \ldots, x_{n})$ be any local coordinate system in $E$. Then, there exists $\{y_{1}, \ldots, y_{r}\} \subset \{x_{1}, \ldots, x_{n}\}$ such that $D_{\alpha}^v = \partial/\partial y_{\alpha}$ are vertical vectors, where $r = \mathrm{dim} VE$. Moreover, $\{D_{1}^v, \ldots, D_{r}^v\}$ span $VE$.
\end{lemma}
\begin{proof}
 Let $(x_{1}, \ldots, x_{n})$ be a local coordinate system in $M$. Let $V$ be a non null vertical vector field on $E$. We can write
\[
 V(x) = a_{i}(x)D_{i}.
\]
Suppose that every $D_{i}$, $i=1, \ldots n$,  is not vertical.  Applying $\pi_{*}$ at $V$ we see that
\[
 \pi_{*}(V(x)) = a_{i}(x)\pi_{*}(D_{i}).
\]
But $\pi_{*}(V(x))=0$, while $a_{i}(x)\pi_{*}(D_{i}) \neq 0$. It is a contradition. Therefore there is a major $r <n$ such that $\{y_{1}, \ldots, y_{r}\}\subset \{x_{1}, \ldots, x_{n}\}$, $D_{\alpha}^v = \frac{\partial}{\partial y^{\alpha}}$, $\alpha=1, \ldots, r$, are vertical and
\[
 V(x) = a_{\alpha}(x)D_{\alpha}^{v}.
\]
It is immediate that $r \leq \mathrm{dim} VE$. If  $r < \mathrm{dim} VE$ we would have a non null vertical vector $U$ such that $U$ is not a linear combination of $\{D_{1}^v, \ldots, D_{r}^v\}$. So $U = a_{j}D_{j}$, where $D_{j} \notin \{D_{1}^v, \ldots, D_{r}^v\}$, $j=1, \ldots, n-r$. As $U$ is vertical we have that some $D_{j}$ would be a  vertical vector. It is a contracdition with the choice of $r$.
\eop
\end{proof}

Let us denote $\tau_{x} \bar{E}$ the set of second order vector field $L$ written as
\[
 L = l^{i\alpha} D_{i\alpha} + l^{\alpha}D_{\alpha},
\]
where $(x_{1}, \ldots, x_{n})$ is a local coordinate system and $(y_{1}, \ldots, y_{r})$ is given by the Lemma above. It is clear that $\tau_{x}\bar{E}$ is a subspace of $\tau_{x}$. Let $\Gamma^{E}$ be a symmetric connection on $M$. We define $\Gamma^{v}: \tau_{x}\bar{E} \rightarrow VE$ by $\Gamma^{v}(L) = \mathbf{v}\Gamma^{E}(L)$. We call $\Gamma^{v}$ of vertical connection on $E$. We observe that $\Gamma^{v}(XV) = \nabla^{v}_{X}V = \mathbf{v}\nabla^{E}_{X}V$, where $X$ and $V$ are vector field and vertical vector field on $E$, respectively.

\begin{definition}
Let $\sigma$ be a section of $\pi$. We define $\bar{Q}_{x}: \tau_{x}\bar{E} \rightarrow TE \odot VE$ by $\bar{Q}(L) = Q(L)$. The section $\alpha^{v}_{\sigma}$ of $\tau^{*}\bar{E} \otimes \sigma^{*}VE$ is given by
\begin{equation}\label{secondformvertical}
 \alpha^{v}_{\sigma} = \Gamma^{v} \mathbf{v}\sigma_{*} - \mathbf{v}\sigma_{*} \Gamma^{M}.
\end{equation}
The vertical second fundamental form of $\sigma$, $\beta_{\sigma}^{v}$ is the unique section of \linebreak $ (TM \odot TM)^*\otimes \sigma^{*}VE$ such that $\alpha_{\sigma}^{v} = \beta^{v}_{\sigma}\circ \bar{Q}$. The tension field of $\beta^{v}_{\sigma}$ is 
\[
 \tau_{\sigma}^{v} = \rm{tr} \beta^{v}_{\sigma}.
\]
We call $\sigma$ a harmonic section if $\tau_{\sigma}^{v}=0$.
\end{definition}

When $\pi$ is a Riemannian submersion, C. Wood, in \cite{wood1}, defines a harmonic section in the same way.

The following linear algebra lemma shows that $\beta_{\sigma}^{v}$ is well defined. 

\begin{lemma}
Let $\alpha^{v}_{\sigma}$ be a section  of $\tau^*E\otimes \sigma^{*}VE$ defined by (\ref{secondformvertical}). Then there exists an unique section $\beta^{v}_{\sigma}$ of $(TM\odot TM)^*\otimes \sigma^{*}VE$ such that $\alpha_{\sigma}^{v}=\beta_{\sigma}^{v} \circ \bar{Q}$.
\end{lemma}
\begin{proof}
Since $\mathrm{Ker}\,\bar{Q} =VE \subset \mathrm{Ker}\,\alpha_{\sigma}^{v}$, the lemma follows from the first isomorphism theorem (see \cite{Rot} pp 67).
\eop
\end{proof}

For the convenience of the reader we repeat the proof of the following Lemma from \cite{catuogno2}.
\begin{lemma}\label{slem1}
Let $\pi:E \rightarrow M$ be an affine submersion with horizontal distribution. Then 
\begin{enumerate}
\item for each $\theta$ in $T^{*}E$, $\int \alpha^{*}_{\sigma} \theta~ d_2X = \frac{1}{2}\int \beta^{*}_{\sigma} \theta (dX,dX);$
\item for each vertical form $\theta$ on $E$, $\int \alpha^{v*}_{\sigma} \theta~ d_2X = \frac{1}{2}\int \beta^{v*}_{\sigma} \theta (dX,dX)$.
\end{enumerate}
\end{lemma}
\begin{proof}
{\it 1.} By definition \ref{formafundamentaleharmonica} of $\beta_{\sigma}$, we have that
\[
\frac{1}{2}\int \beta^{*}_{\sigma} \theta (dX,dX)=\int Q^{*}\beta^{*}_{\sigma} \theta~
d^{2}X = \int (\beta_{\sigma} \circ Q)^{*} \theta~ d^{2}X=\int \alpha^{*}_{\sigma}\theta~ d^{2}X.
\]
The first equality follows from Proposition 6.31 in \cite{emery1}.\\
{\it 2.} The proof is analogous to item 1.
\end{proof}

\begin{proposition}\label{hsprop1}
Let $\sigma$ be a section of $\pi$ and $X$ be a semimartingale in $M$. We have that
\[
 \int \beta_{\sigma}^{*}\theta(dX,dX) = \int \beta_{\sigma}^{v*}\theta(dX,dX)
\]
for all vertical form $\theta$.
\end{proposition}
\begin{proof}
Let $X_{t}$ be a semimartingale in $M$. From Lemma \ref{slem1} we see that
\[
\int \alpha^{*}_{\sigma} \theta d^{2}X=\frac{1}{2} \int \beta_{\sigma}^{*}\theta(dX,dX) \ \ \rm{and} \ \ \int \alpha^{v*}_{\sigma} \theta d^{2}X=\frac{1}{2} \int \beta_{\sigma}^{v*}\theta(dX,dX).
\]
Thus, it is sufficient to show that $\int \alpha^{*}_{\sigma} \theta d^{2}X = \int \alpha^{v*}_{\sigma} \theta d^{2}X$. We first compute $\alpha_{\sigma}(d^{2}X)$. From definition (\ref{fundamental}) we see that
\[
\alpha_{\sigma}(d^{2}X) = \Gamma^E\circ \sigma_{*}(d^{2}X_{t})-\sigma_{*} \circ \Gamma^{M}(d^{2}X_{t}). \\
\]
Write $\sigma_{*} = \mathbf{v}\sigma_{*} + \mathbf{h}\sigma_{*}$. From this, Proposition \ref{pprop1} and Lemma \ref{plem1} we deduce that 
\begin{eqnarray*}
\sigma_{*}(d^{2}X_{t}) 
& = & dX^{i}\sigma_{*}(D_{i}) + d[X^{i},X^{j}]\{\sigma_{*}D_{i},\sigma_{*}D_{j}\}\\
& = & dX^{i}\mathbf{v}\sigma_{*}(D_{i}) + dX^{i}\mathbf{h}\sigma_{*}(D_{i}) + d[X^{i},X^{j}]\{\mathbf{v}\sigma_{*}D_{i},\mathbf{v}\sigma_{*}D_{j}\} \\
& + & 2 d[X^{i},X^{j}]\{\mathbf{v}\sigma_{*}D_{i},\mathbf{h}\sigma_{*}D_{j}\} + 
d[X^{i},X^{j}]\{\mathbf{h}\sigma_{*}D_{i},\mathbf{h}\sigma_{*}D_{j}\}.
\end{eqnarray*}
We observe that $\{\mathbf{v}\sigma_{*}D_{i},\mathbf{v}\sigma_{*}D_{j}\}$ and $ \{\mathbf{v}\sigma_{*}D_{i},\mathbf{h}\sigma_{*}D_{j}\}$ are vertical vectors in $\tau E$, while $\{\mathbf{h}\sigma_{*}D_{i},\mathbf{h}\sigma_{*}D_{j}\}$ are not. Then
\[
\mathbf{v} \sigma_{*}(d^{2}X_{t}) 
 = dX^{i}\mathbf{v}\sigma_{*}(D_{i})+ d[X^{i},X^{j}]\{\mathbf{v}\sigma_{*}D_{i},\mathbf{v}\sigma_{*}D_{j}\} 
 + 2 d[X^{i},X^{j}]\{\mathbf{v}\sigma_{*}D_{i},\mathbf{h}\sigma_{*}D_{j}\}.
\]
Applying the connection $\Gamma^{E}$ at $\sigma_{*}(d^{2}X_{t})$ we see that 
\begin{eqnarray*}
\Gamma^{E} \circ \sigma_{*}(d^{2}X_{t}) 
& = & dX^{i}\mathbf{v}\sigma_{*}(D_{i}) + dX^{i}\mathbf{h}\sigma_{*}(D_{i}) +  d[X^{i},X^{j}]\Gamma^{E}\{\mathbf{v}\sigma_{*}D_{i},\mathbf{v}\sigma_{*}D_{j}\} \\
& + & 2 d[X^{i},X^{j}]\Gamma^{E}\{\mathbf{v}\sigma_{*}D_{i},\mathbf{h}\sigma_{*}D_{j}\} + d[X^{i},X^{j}]\Gamma^{E}\{\mathbf{h}\sigma_{*}D_{i},\mathbf{h}\sigma_{*}D_{j}\}.
\end{eqnarray*}
Since $\mathbf{h}\sigma_{*}D_{i} = H(D_{i})$, from the equality above we conclude that 
\begin{eqnarray*}
\alpha_{\sigma}(d^{2}X_{t})
& = & dX^{i}\mathbf{v}\sigma_{*}(D_{i}) - \mathbf{v} \sigma_{*} \circ \Gamma^{M}(d^{2}X_{t})  + dX^{i}\mathbf{h}\sigma_{*}(D_{i}) - \mathbf{h} \sigma_{*} \circ \Gamma^{M}(d^{2}X_{t})\\
& + & d[X^{i},X^{j}]\Gamma^{E}\{\mathbf{v}\sigma_{*}D_{i},\mathbf{v}\sigma_{*}D_{j}\} + 2 d[X^{i},X^{j}]\Gamma^{E}\{\mathbf{v}\sigma_{*}D_{i},H(D_{j})\} \\
& + & d[X^{i},X^{j}]\Gamma^{E}\{H(D_{i}),H(D_{j})\}. 
\end{eqnarray*}
From assumption (\ref{skewsymmetry}) we see that $\Gamma^{E}\{H(D_{i}),H(D_{j})\}$ are horizontal vectors on $TE$. By definition of $\alpha_{\sigma}^{v}$  we conclude that
\[
 \theta \alpha_{\sigma}(d^{2}X_{t}) = \theta \alpha^{v}_{\sigma}(d^{2}X_{t}).
\]
Consequentely,
\[
 \int \alpha_{\sigma}^{*}\theta d^{2}X = \int \alpha_{\sigma}^{v*}\theta d^{2}X.
\]
Therefore the conclusion follows.
\eop
\end{proof}

\begin{theorem}\label{hsteo1}
 Let $\sigma$ be a section of $\pi$. Then $\sigma$ is a harmonic section if and only if $\mathbf{v}\tau_{\sigma}$ is null.
\end{theorem}
\begin{proof}
 Let $B$ be a Brownian motion in $M$ and $\theta$ be a vertical form in $E$. From Proposition \ref{hsprop1} and the definition of Brownian motion (\ref{Brownian}) we deduce that 
\[
\int \tau_{\sigma}^{*}\theta(B_{t})dt. = \int \tau_{\sigma}^{v*}\theta(B_{t})dt.
\]
Suppose that $\mathbf{v}\tau_{\sigma} =0$. Thus
\[
\int \tau_{\sigma}^{v*}\theta(B_{t})dt.= 0. 
\]
As $\theta$ and $B$ are arbitrary we have $\tau_{\sigma}^{v}=0$. Therefore, by definition, $\sigma$ is a harmonic section. The converse follows from the same argument.
\eop
\end{proof}
\begin{corollary}\label{hscor1}
Let $\sigma$ be a section of $\pi$. If $\sigma$ is a harmonic map, then $\sigma$ is a harmonic section.
\end{corollary}

Let $X,Y$ be vector fields on $M$ and $\sigma$ a section of $\pi$. Then we can write the second fundamental form $\beta_{\sigma}$ and the vertical second fundamental form $\beta_{\sigma}^{v}$ as 
\[
 \beta_{\sigma}(X,Y) = \nabla^{E}_{\sigma_{*}X}\sigma_{*}Y - \sigma_{*}\nabla^{M}_{X}Y
\]
and
 \[
 \beta^{v}_{\sigma}(X,Y) = \nabla^{v}_{\mathbf{v}\sigma_{*}X} \mathbf{v}\sigma_{*}Y - \mathbf{v}\sigma_{*}(\nabla^{M}_{X}Y).
\]

\begin{proposition}\label{hsprop2}
Let $\sigma$ be a section of $\pi$ such that $\sigma$ is not constant. Then the horizontal e vertical components of $\tau_{\sigma}$ are given by
\begin{eqnarray*}\label{hsprop2eq1}
 \mathbf{v}\tau_{\sigma} & = &  \tau_{\sigma}^{v} + \rm{tr}\, C_{\sigma}\\
 \mathbf{h}\tau_{\sigma} & = & \rm{tr}\, D_{\sigma},
\end{eqnarray*}
where 
\begin{equation}\label{hsprop2eq2}
\begin{array}{ccl}
C_{\sigma}(X,Y) & = & 2T_{\mathbf{v}\sigma_{*}X}H(Y) + A_{H(X)}H(Y) + [H(X),\mathbf{v}\sigma_{*}Y]\\
D_{\sigma}(X,Y) & = & T_{\mathbf{v}\sigma_{*}X}\mathbf{v}\sigma_{*}Y + 2A_{H(X)}\mathbf{v}\sigma_{*}Y.
\end{array}
\end{equation}
for  $X,Y$ vector fields on $M$, with $H$ the horizontal lift from $M$ on $E$.
\end{proposition}
\begin{proof}
From definition of $\beta_{\sigma}$ we see that
\[
\beta_{\sigma}(X,Y) = \nabla^{E}_{\sigma_{*}X} \sigma_{*}Y - \sigma_{*}(\nabla^{M}_{X}Y),
\]
for $X,Y$ vector fields on $M$. Write $\sigma = \mathbf{v}\sigma_{*} \oplus \mathbf{h}\sigma_{*}$. We compute the vertical and horizontal components of $\beta_{\sigma}$. First, we develop the vertical component. Thus
\begin{eqnarray*}
\mathbf{v}\beta_{\sigma}(X,Y) & = & \mathbf{v}\nabla^{E}_{\mathbf{v}\sigma_{*}X} \mathbf{v}\sigma_{*}Y  +\mathbf{v}\nabla^{E}_{\mathbf{v}\sigma_{*}X} \mathbf{h}\sigma_{*}Y +\mathbf{v}\nabla^{E}_{\mathbf{h}\sigma_{*}X} \mathbf{v}\sigma_{*}Y + \mathbf{v}\nabla^{E}_{\mathbf{h}\sigma_{*}X} \mathbf{h}\sigma_{*}Y \\
& - & \mathbf{v}\sigma_{*}(\nabla^{M}_{X}Y).
\end{eqnarray*}
From the definitions of the Fundamental tensors $T$ and $A$ we deduce that
\begin{eqnarray*}
\beta_{\sigma}(X,Y)  
& = & \beta_{\sigma}^{v}(X,Y) + 
2T_{\mathbf{v}\sigma_{*}Y} \mathbf{h}\sigma_{*}X + A_{\mathbf{h}\sigma_{*}X} \mathbf{h}\sigma_{*}Y - \mathbf{v}[\mathbf{h}\sigma_{*}X,\mathbf{v}\sigma_{*}Y].
\end{eqnarray*}
Let $\{e_{1}, \ldots, e_{n}\}$ be a local orthonormal system at $x \in M$ and $H(e_{i})$, $i=1, \ldots, n$, their horizontal lifts on $E$. It is clear that $H(e_{i})=\mathbf{h}\sigma_{*}e_{i}$. A direct computation shows that
\[
\mathrm{tr}\,\mathbf{v}\beta_{\sigma}  
= \mathrm{tr}\,\beta^{v}_{\sigma} + \sum_{i=1}^{n}(2T_{\mathbf{v}\sigma_{*}e_{i}}H(e_{i}) + A_{H(e_{i})}H(e_{i}) - [H(e_{i}),\mathbf{v}\sigma_{*}e_{i}]).
\]
Writting  
\[
C_{\sigma}(X,Y) = 2T_{\mathbf{v}\sigma_{*}X}H(Y) + A_{H(X)}H(Y) + [H(X),\mathbf{v}\sigma_{*}Y],
\]
where $X,Y$ are vector field in $M$, we conclude that 
\[
 \mathbf{v}\tau_{\sigma} = \tau^{v}_{\sigma} +\rm{tr}\, C_{\sigma}.
\]
Using the fact that $\pi$ is an affine submersion with horizontal distribution it is possible to show that 
\[
 \mathbf{h}\tau_{\sigma} = \rm{tr}\, D_{\sigma},
\]
where
\[
D_{\sigma}(X,Y) = T_{\mathbf{v}\sigma_{*}X}\mathbf{v}\sigma_{*}Y + 2A_{H(X)}\mathbf{v}\sigma_{*}Y.
\]
\eop
\end{proof}

\begin{theorem}\label{hsteo2}
Let $\sigma$ be a section of $\pi$ and $C_{\sigma}$ and $D_{\sigma}$ be tensors defined by (\ref{hsprop2eq2}). 
\begin{enumerate}
\item If $\rm{tr}\, D_{\sigma}= 0$, then $\sigma$ is a harmonic map if and only if $\sigma$ is a harmonic section;
\item If $\sigma$ is a harmonic map, then $\rm{tr}\, C_{\sigma} = 0 $ and $\rm{tr}\, D_{\sigma}=0$;
\item If $\rm{tr}\, C_{\sigma} \neq 0$ or $\rm{tr}\, D_{\sigma} \neq 0$,  then $\sigma$ is not a harmonic map;
\item If $\sigma$ is a harmonic section, then $\rm{tr}\, C_{\sigma} = 0$;
\item If $\rm{tr}\, C_{\sigma} \neq 0$, then $\sigma$ is not a harmonic section.
\end{enumerate}
\end{theorem}
\begin{proof}
{\it 1.} It follows immediately from Theorem \ref{hsteo1} and Proposition \ref{hsprop2}.\\ 
{\it 2.} Suppose that $\sigma$ is a harmonic map. Then $\mathbf{v}\tau_{\sigma} = 0$ and $\mathbf{h}\tau_{\sigma} =0$. From Theorem \ref{hsteo1} we have that $\tau^{v}_{\sigma}=0$. Proposition \ref{hsprop2} now shows that $\rm{tr}\,C_{\sigma} =0$ and $\rm{tr}\, D_{\sigma}=0$.\\
{\it 3.} If $\rm{tr}\, D_{\sigma} \neq 0$, the result follows immediatelly. Suppose that $\rm{tr}\, C_{\sigma} \neq 0$. If $\sigma$ were a harmonic map, we would have $\mathbf{v}\sigma_{*} = 0$. Theorem \ref{hsteo1} would show that $\tau^{v}_{\sigma}=0$. From Proposition \ref{hsprop1} we would conclude that $\rm{tr}\, C_{\sigma} = 0$. It is a contracdition. \\
{\it 4.} Suppose that $\sigma$ is a harmonic section. By definition $\tau^{v}_{\sigma}=0$. Therefore, from Theorem \ref{hsteo1} we have that $\mathbf{v}\tau_{\sigma} = 0$. Proposition \ref{hsprop2} now shows that $\rm{tr}\, C_{\sigma}=0$.\\
{\it 5.} We assume that $\rm{tr}\,C_{\sigma} \neq 0$. Suppose that $\sigma$ is a harmonic section. Theorem \ref{hsteo1} shows that $\mathbf{v}\tau_{\sigma}=0$. As $\tau^{v}_{\sigma} =0$ we have $\rm{tr}\, C_{\sigma} = 0$. It is a contracdition.
\eop\\
\end{proof}

\section{Examples}

{\large {\it Product manifold}}\\

Let $M$ be a Riemannian manifold, $N$ a differential manifold and $E = M \times N$ a product manifold. Let $\nabla^{M}$ be the Levi-Civita connection on $E$ and $\nabla^{N}$ a symmetric connection on $N$. Furthermore, we endow $E$ with the product connection $\nabla^{E} = \nabla^{M} + \nabla^{N}$. We observe that $\nabla^{E}$ is also symmetric. It follows immediately that $\pi: E \rightarrow M$ is an affine submersion with horizontal distribution. It is also easy to check that $T \equiv 0$ and $A \equiv 0$. From this and Theorem \ref{hsteo2} we deduce the following.
\begin{proposition}
 Let $\pi:E=M\times N \rightarrow M$ be a projection. Let $\sigma$ be a section of $\pi$. Then $\sigma$ is a harmonic map if and only if it is a harmonic section.
\end{proposition}
Furthermore, it follows easily that $\rm{tr}\, C_{\sigma} =0$. \\

{\large {\it Tangent Bundle with horizontal Lift}}\\

Let $M$ be a Riemmanian manifold and $TM$ its tangent bundle. Let $\nabla^{M}$ be a symmetric connection on $M$. It is possible to prolong $\nabla^{M}$ to a connection on $TM$. A well known way is the horizontal lift $\nabla^{h}$ (see \cite{yano} for the definition of $\nabla^{h}$). We observe that $\nabla^{M}$ is symmetric and flat if and only if $\nabla^{h}$ is symmetric (see proposition 7.3, chapter II in \cite{yano}). Let $X,Y$ be vector fields on $M$, so $\nabla^{h}$ satisfies the following equations:
\begin{equation}
\begin{array}{ccl} 
\displaystyle \nabla^{h}_{X^{V}}Y^{V} & = & 0  \\
\nabla^{h}_{X^{V}}Y^{H} & = & 0  \\
\nabla^{h}_{X^{H}}Y^{V} & = & (\nabla_{X}Y)^{V} \\
\nabla^{h}_{X^{H}}Y^{H} & = & (\nabla_{X}Y)^{H},
\end{array}
\end{equation}
where $X^{V}, Y^{V}$ are vertical lifts, while $X^{H},Y^{H}$ are horizontal lifts (see \cite{yano} for more details). From this we see that $ T \equiv 0$ and $A \equiv 0$. It follows that $\rm{tr}\,D_{\sigma} = 0$. So we have showed the following.
\begin{proposition}
 Let $M$ be a Riemmanian manifold endowed with a flat Levi-Civita connection $\nabla^{M}$. Let $TM$ be its tangente bundle endowed with horizontal lift $\nabla^{h}$. Let $\sigma$ be a section of $\pi:TM \rightarrow M$. Then $\sigma$ is harmonic map if and only if $\sigma$ is a harmonic section.
\end{proposition}

Moreover, from the fact that $T \equiv 0$ and $A \equiv 0$ we see that
\[
\rm{tr}\,C_{\sigma} = - [e_{i}^{H},\mathbf{v}\sigma_{*}(e_{i})].\\
\]

{\large {\it The affine submersion in the sense of Blumenthal}}\\

 Let $M$ be a Riemmanian manifold and $E$ a differential manifold. We endow $M$ and $E$ with symmetric connections $\nabla^{M}$ and $\nabla^{E}$. Blumenthal called a submersion $\pi:(E,\nabla^{E}) \rightarrow (M, \nabla^{M})$ an affine submersion if $\pi_{*}$ commutes with parallel translations induced by affine connections on $E$ and $M$, respectively. N. Abe and K. Hasegawa  showed, see Theorem 5.1 in \cite{abe}, that $\pi$ is an affine submersion in the sense of Blumenthal if and only if for an arbitrary horizontal bundle, $\pi:(E, \nabla^{E}) \rightarrow (M, \nabla^{M})$ is an affine submersion such that $T_{V}W = 0$, $A_{Y}W = 0$ and $\mathbf{h}\nabla^{E}_{V}H(X) =0$, where $V,W$ are vertical vector fields on $E$, $Y$ is a horizontal vector field on $E$ and $X$ is a vector field on $M$. From this and Theorem \ref{hsteo2} we conclude the following.
\begin{proposition}
 if $\pi:(E,\nabla^{E}) \rightarrow (M, g)$ is a submersion in the sense of Blumenthal, where $\nabla^{E}$ is a symmetric connection and $M$ is a Riemmanian manifold, then a section $\sigma$ of $\pi$ is a harmonic map if and only if $\sigma$ is a harmonic section.
\end{proposition}

{\large {\it Riemannian submersion}}\\

 Let $\pi:E \rightarrow M$ be a Riemannian submmersion with totally geodesic fibers. We stated some results about Fundamental tensors $T$ and $A$ in this case. Let $X,Y$ be horzontal fields on $E$ and $U,V$ be vertical fields on $E$. B. O'Neill showed in \cite{oneill} that 
\[
(i)\, (T_{U}V,X) = - (T_{U}X,V) \ \ (ii)\, A_{X}Y = \frac{1}{2}\mathbf{v}[X,Y] \ \ (iii)\, (A_{X}Y,U) = - (A_{X}U,Y)
\]
From $(i)$  we see that the totally geodesic fibers property is equivalent to $T \equiv 0$. From this and $(ii)$ we conclude that $\rm{tr}\, C_{\sigma} = \sum_{i=0}^{n}[H(e_{i}), \mathbf{v}\sigma_{*}(e_{i})]$ and $\rm{tr}\, D_{\sigma} = \sum_{i=1}^{n} 2A_{H(e_{i})}\mathbf{v}\sigma_{*}(e_{i})$ for any local orthonormal system $\{e_{1}.\ldots, e_{n}\}$ on $M$. As $H$ is an isomorphism from $TM$ into $HE$ we have that $\{H(e_{1}), \ldots, H(e_{n})\}$ are basis in $HE$. Therefore any horizontal field $Y$ on $E$ can be written as $Y = a_{j} H(e_{j})$. From $(ii)$ and $(iii)$ we deduce that
\begin{eqnarray*}
 (\sum_{i=1}^{n}A_{H(e_{i})}\mathbf{v}\sigma_{*}(e_{i}), Y) & = & -\sum_{i=1}^{n}(A_{H(e_{i})}Y,\mathbf{v}\sigma_{*}(e_{i}))\\
& = & -\sum_{i,j=1}^{n}\frac{1}{2}a_{j}(\mathbf{v}[H(e_{i}),H(e_{j})],\mathbf{v}\sigma_{*}(e_{i})). 
\end{eqnarray*}

Thus, if $\mathbf{v}[H(e_{i}),H(e_{j})]$ is null for all $i,j=1,\ldots,n$, then  $\sum_{i=1}^{n}A_{H(e_{i})}\mathbf{v}\sigma_{*}(e_{i}) = 0$. As $\mathbf{v}[\cdot,\cdot]$ is a tensor we have the following.

\begin{proposition}\label{exprop4}
 Let $\pi:E \rightarrow M$ be a Riemannian submmersion with totally geodesic fibers. Suppose that the horizontal distribution given by $HE$ is integrable. Let $\sigma$ be a section of $\pi$. Then $\sigma$ is a harmonic map if and only if it is a harmonic section.
\end{proposition}
Furthermore, from the fact that $T\equiv 0$ and $(ii)$ we conclude that 
\begin{equation}\label{exeq1}
\rm{tr}\, C_{\sigma} = - \sum_{i=0}^{n}[H(e_{i}), \mathbf{v}\sigma_{*}(e_{i})].\\
\end{equation}

\end{document}